# Lyapunov-based Model Reference Adaptive Controller Design for a Class of Nonlinear Fractional Order Systems

Seyed Mohammad Moein Mousavi
Student
Electrical and Computer Engineering Department
Tarbiat Modares University
Tehran, Iran
Email: moein_mousavi@modares.ac.ir

Mohammad T.H. Beheshti
Professor
Electrical and Computer Engineering Department
Tarbiat Modares University
Tehran, Iran
Email: mbehesht@modares.ac.ir

Amin Ramezani
Professor
Electrical and Computer Engineering Department
Tarbiat Modares University
Tehran, Iran
Email: ramezani@modares.ac.ir

*Abstract*— **This paper is concerned with model reference adaptive controller design for a class of nonlinear fractional order systems. Recent works on this topic rarely include direct methods and they are mostly based on indirect methods where the frequency distributed model is used to prove the stability of the closed loop system. Since the chain rule cannot be applied in fractional derivations, in order to prove the lyapunov stability here fractional inequalities are used. Finally, by means of a numerical example, the controller performance is demonstrated.**

## 1 Introduction

Despite the long history of fractional calculus, it is recently being utilized in control engineering which is mainly because of its special features for better modelling some physical systems. There are some systems like a thermal system, which can be described better using fractional order models compare to rational order ones [8]. Besides, as the real systems are mostly time variant, using adaptive controllers is important because these controllers are robust to parameters variations during the time.

Hence, adaptive controller design for linear and nonlinear fractional order models is a significant subject which is being considered in the very recent years.in [4] model reference adaptive controller is designed using indirect method.in [3], using grunwald definition of fractional order derivatives, fractional order model is converted to an integer order discrete time model and the adaptive controller is designed using STR technique. In [2,9] a new method for parameter estimation in nonlinear fractional order systems is introduced and the controller parameters are computed using system parameters (indirect method). Then the stability of tracking and estimation error is proved using frequency distributed model. In [7,11,12], adaptive laws for a 1-dimensional nonlinear fractional order systems is introduced and the stability of tracking error is proved using frequency distributed model. but this work is not completed for n-dimension case of nonlinear fractional order systems. besides, lyapunov based controller design (obtaining the adaptation laws using direct lyapunov method) is considered only for the linear case in [6,10].

To the best of our knowledge, lyapunov based model reference adaptive controller design for a class of nonlinear n-dimensional fractional order systems is an open problem which is not discussed yet.in the rest of this paper in section 2, some preliminaries on fractional calculus and different definitions of fractional derivatives are introduced.in section 3 , lyapunov based controller is designed and the stability is proved.in section 4 a numerical simulation is introduced to depict the effectiveness of the proposed controller and in section 5 we have the conclusion part.

## 2 Preliminaries

Let $C[a\ b]$ be the space of continuous functions $f(t)$ on $[a\ b]$ and we mean by $C^k$ the space of real-valued functions $f(t)$ with continuous derivatives up to order $k-1$ such that $f^{(k-1)}(t) \in C[a\ b]$ and $f^i(t)$ is the i-th derivative of $f(t)$.

### 2.1 fractional calculus

According to [1] there are three main definitions of fractional order derivatives:

- The Riemann–Liouville fractional derivative of order $\alpha$ of $f(t) \in C^m[a\ b]; t \in [a\ b]$:

$$^{RL}_a D_t^\alpha f(t) = \frac{1}{\Gamma(m-\alpha)} \frac{d^m}{dt^m} \times \int_a^t (t-\tau)^{m-\alpha-1} f(\tau) d\tau \qquad (1)$$

- Caputo's derivative of order $\alpha$ of $f(t) \in C^m[a\ b]; t \in [a\ b]$:

$$^C_a D_t^\alpha f(t) = \frac{1}{\Gamma(m-\alpha)} \int_a^t (t-\tau)^{m-\alpha-1} \frac{d^m f(\tau)}{dt^m} d\tau \qquad (2)$$

- Grunwald-Letnikov definition:

$$^{GL}_aD^\alpha_t f(t) = \lim_{h\to\infty} \frac{1}{\Gamma(a)h^\alpha} \sum_{j=0}^{[\frac{t-a}{h}]} \frac{\Gamma(a+j)}{\Gamma(j+1)} f(t-jh) \quad (3)$$

Where $\left[\frac{t-a}{h}\right]$ denotes the integer part of $\frac{t-a}{h}$ and $\Gamma(\ )$ is Euler gamma function.

The drawback of the first definition is that the initial conditions are in terms of the variable's fractional order derivatives. having said that, Caputo's definition of fractional derivative needs the initial conditions of the main function and not its fractional derivatives. hence, in engineering usages caputo definition is commonly applied. for the simplicity in the rest of this paper this notation is used:

$$^C_0D^\alpha_t f(t) \triangleq D^\alpha f(t) \quad (4)$$

*2.2  fractional inequalities*

These two lemmas are proven in [5] which will help us with stability analyze in fractional order systems:

**Lemma1**: Let $x(t) \in \mathbb{R}^n$ be a differentiable vector. For any $t \geq t_0$:

$$\frac{1}{2} D^\alpha(e^T P e) \leq e^T P D^\alpha e \quad (5)$$

Where $P \in^{n \times n}$ is a symmetric positive definite and constant matrix

**Lemma 2:** For any $t \geq t_0$

$$\frac{1}{2} D^\alpha(tr(A^T A)) \leq 2tr(A^T D^\alpha A), \forall \alpha \in (0,1] \quad (6)$$

## 3  CONTROLLER DESIGN

Let us consider the nonlinear fractional order systems as:

$$D^\alpha x(t) = Ax(t) + Bu(t) + F(x(t)) \quad (7)$$

Where $A \in \mathbb{R}^{n \times n}$ and $B \in \mathbb{R}^{n \times 1}$ are unknown constant matrixes and the pair $[A\ B]$ is controllable. $x$ is a n-dimensional state, $u \in \mathbb{R}^{1 \times 1}$ is the control input and:

$$F(x(t)) = [f_1(x) \quad f_2(x) \quad \cdots \quad f_n(x)]^T \quad (8)$$

Where $f_i(x)$ is presumed to be a Lipschitz functions.

And the desired model is:

$$D^\alpha x_m(t) = A_m x_m(t) + B_m r(t) \quad (9)$$

Where $A_m \in \mathbb{R}^{n \times n}$ and $B_m \in \mathbb{R}^{n \times 1}$ are constant matrixes, $A_m$ is Hurwitz, $r(t) \in \mathbb{R}^{1 \times 1}$ is the reference signal, And the pear $[A_m, B_m]$ is controllable Our objective here is to force the systems (7) to follow the response of system (9) so that all the state errors converge to zero in finite time. We consider the Input signal as:

$$u(t) = \theta_1 x(t) + \theta_2 r(t) + \theta_3 F(x(t)) \quad (10)$$

Hence we have:

$$D^\alpha x(t) = (A + B\theta_1)x + B\theta_2 r + (B\theta_3 + I)F \quad (11)$$

And the error is defined as:

$$e = x - x_m \quad (12)$$

So the error dynamic becomes:

$$D^\alpha e = A_m e + (A + B\theta_1 - A_m)x + (B\theta_2 - B_m)r + (B\theta_3 + I)F \quad (13)$$

Nominal values of the controller gains are taken as $\theta_{1_0}, \theta_{2_0}, \theta_{3_0}$ and we define:

$$J = \theta_1 - \theta_{1_0} \quad (14)$$
$$\theta_{1_0} = B^{-1}(A_m - A)$$

$$K = \theta_2 - \theta_{2_0} \quad (15)$$
$$\theta_{2_0} = B^{-1} B_m$$

$$L = \theta_3 - \theta_{3_0} \quad (16)$$
$$\theta_{3_0} = -B^{-1}$$

So we have:

$$D^\alpha e = A_m e + BJx + BKr + BLF \quad (17)$$

Now a Lyapunov candidate is adopted as:

$$V = \frac{1}{2}[e^T P e + tr(J^T J + K^T K + L^T L)] \quad (18)$$

Which P is a symmetric positive definite matrix. Using lemmas 1 and 2 We have:

$$D^\alpha V \leq e^T P D^\alpha e + tr(J^T D^\alpha J + K^T D^\alpha K + L^T D^\alpha L) = e^T P A_m e + e^T P B J x + e^T P B K r + e^T P B L F + tr(J^T D^\alpha J + K^T D^\alpha K + L^T D^\alpha L) \quad (19)$$

Since $e^T P A_m e$ is scalar, we can write:

$$e^T P A_m e = \frac{1}{2}(e^T P A_m e + e^T A_m^T P e) = \frac{1}{2} e^T (P A_m + A_m^T P) e \quad (20)$$

And Since $A_m$ is horowitz, there exist a positive definite $Q$ such that:

$$P A_m + A_m^T P = -Q \quad (21)$$

So if these equations are valid, the error dynamic would be stable:

$$e^T P B J x = -tr(J^T D^\alpha J) \quad (22)$$

$$e^T P B K r = -tr(K^T D^\alpha K) \quad (23)$$

$$e^T P B L F = -tr(L^T D^\alpha L) \quad (24)$$

$e^T P B J x$  Is scalar. so we can write:

Similarly:
$$e^T PBKr = tr(re^T PBK) \quad (26)$$
$$e^T PBLF = tr(xe^T PBL) \quad (27)$$

So with these adaptation laws, the error dynamic would be stable:
$$\begin{cases} D^\alpha \theta_1 = -B^T P e x^T \\ D^\alpha \theta_2 = -B^T P e r^T \\ D^\alpha \theta_3 = -B^T P e F^T \end{cases} \quad (28)$$

Since the unknown parameters cannot appear in the adaptation laws, like 1-dimensional case where the unknown parameter is merged with the adaptation gain, here we can also replace $B$ with $B_m$ which is known, with assumption that $B$ and $B_m$ have the same signs.

## 4 NUMERICAL SIMULATIONS

Consider this unstable nonlinear fractional order system:
$$\begin{cases} D^{0.7} x_1 = x_2 \\ D^{0.7} x_2 = x_1 + x_2 + x_1^2 + \sin(x_2) + u \end{cases}$$

It is assumed that this systems' parameters are unknown but the sign of matrix $B$ in (7) is known. And The desired model is:
$$\begin{cases} D^{0.7} x_{1m} = x_{2m} \\ D^{0.7} x_{2m} = -5x_{1m} - 5x_{2m} + 5r \end{cases}$$

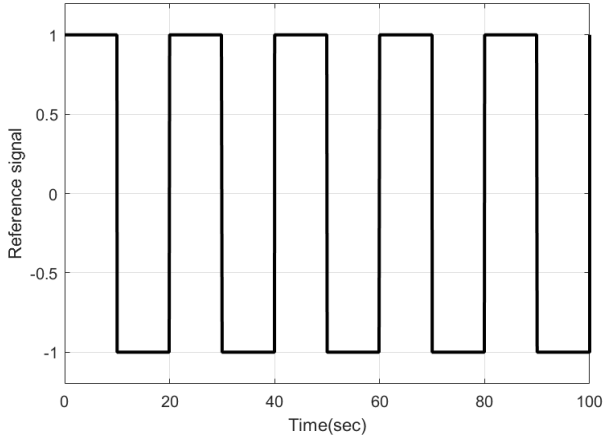

Figure 1. The reference signal

If we take:
$$P = \begin{pmatrix} 20 & 10 \\ 10 & 20 \end{pmatrix} > 0$$

Then:
$$Q = \begin{pmatrix} 100 & 130 \\ 130 & 180 \end{pmatrix} > 0$$

The reference signal is a square wave with period 20 sec which is shown in Figure 1. Simulation is performed using FOMCON toolbox [13]. If we apply the obtained adaptation laws (28) then following figures are derived.

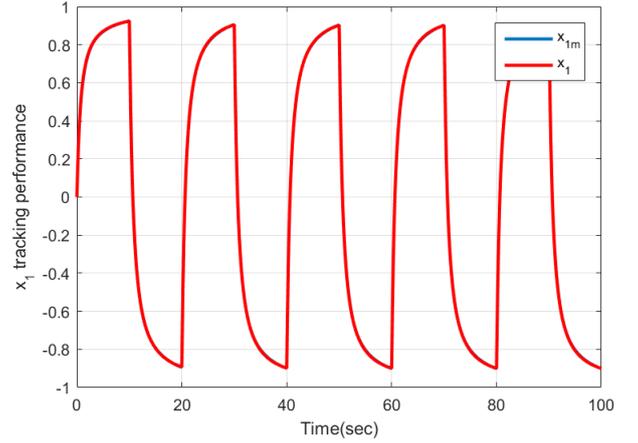

Figure 2. tracking the first state

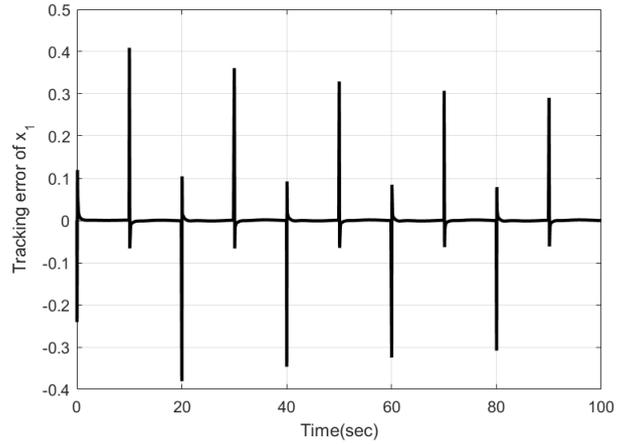

Figure 3. $x_1$ tracking error

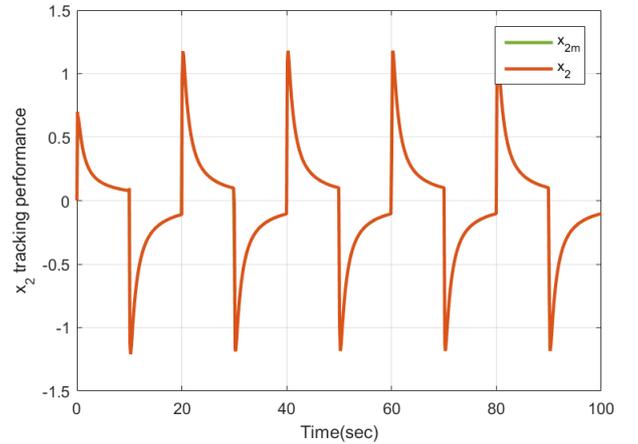

Figure 4. tracking the second state

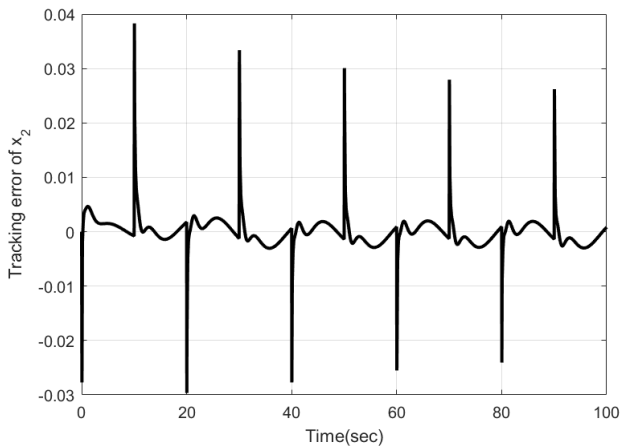

Figure 5. $x_2$ tracking error

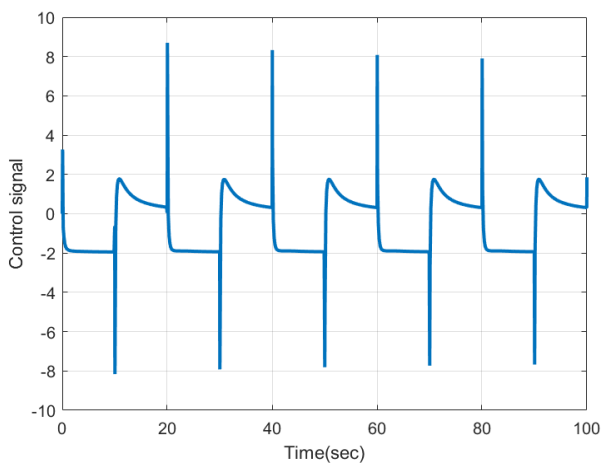

Figure 6. Control signal

Figures 2,3 depict the first state of the desired and controlled model and the tracking error. Despite this state is the output in both systems, from section 3 we expect that all the state errors converge to zero in finite time and Figures 4,5 represent this fact. So One can note that both errors converge to zero asymptotically and the two states are controlled by means of designed controller. Also we can see control signal in Figure 6.

## 5  CONCLUSION

In this paper an adaptive controller was designed for a class of nonlinear fractional order systems based on direct lyapunov method and for this aim, we used fractional inequalities. In a numerical example an unstable nonlinear fractional order system was forced to follow a stable linear fractional order system with satisfactory time constant and through this example we illustrated the effectiveness of the proposed controller. A future study can be:

- Controller design in a way which does not depend on nonlinear dynamic of the system
- None zero initial conditions case

REFERENCES

[1] Podlubny I. Fractional differential equations: an introduction to fractional derivatives, fractional differential equations, to methods of their solution and some of their applications. Academic press; 1999.

[2] Y.Chen,Y.Wei,S.Liang,Y.Wang, "Indirect model reference adaptive control for a class of fractional order systems", Communications in Nonlinear Science and Numerical Simulation Volume 39, October 2016, Pages 458–471

[3] M.Abedini,M.Nojoumian,H.Salarieh,A.Meghdari,"Model reference adaptive control in fractional order systems using discrete time approximation methods",Communications in Nonlinear Science and Numerical Simulation, Volume 25, Issues 1–3, August 2015, Pages 27–40

[4] S.Ladaci,Y.Bensafia, "Indirect fractional order pole assignment based adaptive control", Engineering Science and Technology, an International Journal Volume 19, Issue 1, March 2016, Pages 518–530

[5] N.Aguila-Camacho,M.A.Duarte-Mermoud,J.A.Gallegos,"Lyapunov functions for fractional order systems", Communications in Nonlinear Science and Numerical Simulation Volume 19, Issue 9, September 2014, Pages 2951–2957

[6] M.A.Duarte-Mermoud,N.Aguila-Camacho,J.A.Gallegos,R.Castro Linares,"Using general quadratic Lyapunov functions to prove Lyapunov uniform stability for fractional order systems", Communications in Nonlinear Science and Numerical SimulationVolume 22, Issues 1–3, May 2015, Pages 650–659

[7] B.Shi,J.Yuan,C.Dong,"On Fractional Model Reference Adaptive Control", The Scientific World Journal Volume 2014 (2014), Article ID 521625, 8 pages

[8] A.Aribi, C.Farges, M.Aoun, P.Melchior, S.Najar, M.Abdelkrim, "Fault detection based on fractional order models: Application to diagnosis of thermal systems", Communications in Nonlinear Science and Numerical Simulation Volume 19, Issue 10, October 2014, Pages 3679–3693.

[9] Y.Chen,S.Cheng,Y.Wei,Y.Wang,"Indirect model reference adaptive control for a class of linear fractional order systems", 2016 American Control Conference (ACC) Boston Marriott Copley Place July 6-8, 2016. Boston, MA, USA

[10] H.Cong,Q.Zhidong,M.Qian,Z.Xi,"Factional Order Model Reference Adaptive Control based on Lyapunov Stability Theory", Proceedings of the 35th Chinese Control Conference July 27-29, 2016, Chengdu, China

[11] Y.Wei,S.Liang,Y.Hu,Y.Wang," Composite model reference adaptive control for a class of nonlinear fractional order systems",Proceedings of the ASME 2015 International Design Engineering Technical Conferences & Computers and Information in Engineering Conference IDETC/CIE 2015 August 2-5, 2015, Boston, Massachusetts, USA

[12] Y.Wei,Y.Hu,L.Song,Y.Wang, "Tracking Differentiator Based Fractional Order Model Reference Adaptive Control: The $1 < \alpha < 2$ Case" 53rd IEEE Conference on Decision and Control December 15-17, 2014. Los Angeles, California, US

[13] A. Tepljakov, E. Petlenkov, and J. Belikov, "FOMCON: a MATLAB toolbox for fractional-order system identification and control," International Journal of Microelectronics and Computer Science, vol. 2, no. 2, pp. 51–62, 2011.